\documentclass[12pt, amsmath, amsfonts]{article}

\usepackage{amscd}
\usepackage{amssymb}
\def\K{\mathop {\mathrm {K}}\nolimits}
\def\KK{\mathop {\mathrm {KK}}\nolimits}

\def\Mat{\mathop{\mathrm {Mat}}\nolimits}

\newtheorem{thm}{Theorem}[section]
\newtheorem{defn}[thm]{Definition}
\newtheorem{prop}[thm]{Proposition}
\newtheorem{lem}[thm]{Lemma}
\newtheorem{cor}[thm]{Corollary}
\newtheorem{rem}[thm]{Remark}

\def\proof{{\sc Proof}.\ }
\def\endproof{\hfill$\square$}

\em
\begin{document}
\title{Noncommutative Spherical Tight Frames in finitely generated Hilbert C*-modules\footnote{The work was supported in part by Vietnam National Project of Research in  Fundamental Sciences and The Abdus Salam ICTP, UNESCO.}}
\author{Do Ngoc Diep}
\maketitle
\begin{abstract} Let $A$ be a fixed C*-algebra.
In an arbitrary finitely generated projective $A$-module $V\subseteq A^n$, a spherical tight $A$-frame is a set of of $k$, $k>n$, elements $f_1,\dots,f_k$ such that the associated matrix $F = [f_1, \dots,f_k]$ up-to a constant multiple is a partial isometry of the Hilbert structure on the projective finitely generated $A$-module $V$. The space $\mathcal F^A_{k,n}$ of all such $A$-frames form a C*-algebra, generated by a system of partial isometries and the structure of such C*-algebras are well described, especially in the case $A=\mathbf R$ or $\mathbf C$: The main result of K. Dykema and N. Strawn for these cases are generalized to our general projective finitely generated Hilbert $A$-module case. This generalization gives the possibility to study the universal classifying space.
\end{abstract}
\section{Introduction}

In the classical (commutative) theory of vector bundles the construction of Stiefel bundles is well-known and gives rise to the classification space for principal bundles. Let us recall (\cite{husemoller}, Chapter 7) that for ground field $K= \mathbf R, \mathbf C$ or $\mathbf H$, denote $U_K(k)$ the orthogonal group $O_K(k)$ for $K= \mathbf R$, the unitary group $U(k)$ for $K=\mathbf C$ and the symplectic groups $Sp(k)$ for $K=\mathbf H$. 
Denote $$K^1 \subset K^2 \subset \dots \subset K^k \subset K^{k+1} \subset \dots K^\infty := \bigcup_{k=1}^\infty K^k$$ and also $U_K(\infty) = \bigcup_{k=1}^\infty U_K(k)$.  The {\it Stiefel variety}, by definition is the subspace $V_k(K^n)$ of $K^{kn}$ consisting of orthonormal frames ($k$-tuples $[v_1,\dots, v_k]$ of orthonormal vectors) in $K^k$, for $k\leq n = 1,2,\dots,\infty$. One defines the action of $U_K(k)$ on $V_k(K^n)$ by the map $\eta: U_K(k) \to V_k(K^n)$ via the formula $\eta(u) = (u(\mathbf e_1), \dots, u(\mathbf e_k)),$ where $\mathbf e_1,\dots, \mathbf e_k$ is the standard orthonormal basis in $K^n$. With this action of $U_K(k)$  {\it the space $V_k(K^n)$ becomes a principal bundle with structural group $U_K(k)$}. It is well-known that $V_k(K^n) \cong U_K(n)/U_K(n-k)$, in particular it is a homogeneous space, called {\it Stiefel variety} and for example for $k=n$, $V_n(K^n) = U_K(n)$, for $k=1$ $V_1(K^n)$ is the sphere $\mathbf S^m$, where $m= n-1$ 
in the case 
$F= \mathbf R$, $m= 2n-1$ in the case $F= \mathbf C$ and $m= 4n-1$ in the case $F= \mathbf H$, for $k=n-1$, $V_{n-1}(K^n) \cong U_K(n)/U_K(1) \cong SU_K(n).$  

Let us recall also that the Grassman manifold $G_k(K^n)$ is the the space of $k$-dimensional subspaces in $K^n$. It is also well-known, see e.g. (\cite{husemoller}, Theorem 2.2) that $U_K(k) \hookrightarrow V_k(K^n) \twoheadrightarrow G_k(K^n)\cong U_K(n)/U_K(n-k)$ is a principal bundle. Moreover, see e.g. (\cite{husemoller}, Theorem 4.1) the natural inclusion $U_K(n) \to U_K(n+q)$ induces morphism of homotopy groups $\pi_i(U_K(n)) \to \pi_i(U_K(n+q))$ which are isomorphism for $i\leq c(n+1)-3$ and epimorphism for $i\leq c(n+1)-2$. For $q=\infty$ the hypotheses are satisfied. 

From these one deduces, see e.g. (\cite{husemoller}, Chapter 7, Theorem 6.1) that the principal bundle $V_k(K^n) \to G_k(K^n)$ is universal in dimension less than or equal to $c(m+1) -2$ and $V_k(K^\infty) \to G_k(K^\infty)$ is  a universal bundle. This result can be understood, for example in the sense that the space $[X,G_k(K^n)] \cong Vect_k(X) $ of  homotopy classes of maps is isomorphic with the space of isomorphic classes of vector bundles $f^*(\gamma_k^{k+m})$ on a CW complex $X$ obtained from a universal vector bundle $\gamma_k^{k+m}$ on $G_k(K^n)$ when $n \leq c(m+1) -2$. In other words, all vector bundles on a CW complex $X$ can be regarded as some induced vector bundle,  associated with a map from $X$ to the classification space $G_k(K^{k+m})$, for $n\leq c(m+1)-2$.

In the work \cite{dykemastrawn}, K. Dykema and N. Strawn had proved that the space of frames with $k>n$ and with replacement of the orthonormality of the system of vectors $f_i, i\in I \subset \{1,\dots,k\}$ by the spherical tight condition of frames $F= \{f_i\}_{i\in I}$,
$$b||v||^2 = \sum_{i\in I} |\langle v,f_i\rangle |^2,\quad b=const,$$ admits also a manifold structure. The aim of this paper is to use some noncommutative analog of this manifold structure to construct noncommutative universal classifying spaces. 
We  introduce, in Section 2  a noncommutative setting necessary for C*-algebras Hilbert projective $A$-modules of finite type. In Section 3, we proved that the space of all noncommutative spherical tight $A$-frames indeed admits a structure of a C*-algebra, generated by a system of partial isometries and which can be decomposed into a orthogonal sum of simple ones. In Section 4 we apply this result to construct the universal classifying spaces in general noncommutative situation.

\section{Equivalence classes of noncommutative spherical tight $\mathbf A$-frames}

We are interested in a class of noncommutative Serre Fibrations (NCSF) \cite{diep1} of type 
$$\CD f: A @>>> B\endCD,$$
i. e. a morphism from A to B in the category of C*-algebras with NCCW structure. In that case, $B$ is endowed with a structure of $A$-module:
$$a\in A, b\in B \mapsto a.b:= f(a)b.$$
A {\it section} of the noncommutative Serre fibration $\CD f: A @>>> B\endCD,$ is defined as a morphism $s: B \to A$ such that $f\circ s = Id_B$.

\begin{defn}{\rm
A {\it noncommutative spherical tight $A$-frame} is a collection $F=[f_i]_{i\in I}$ of vectors $f_i$ in a finitely generated Hilbert $A$-module satisfying the condition
$$b ||v||_A^2 = \sum_{i\in I} ||\langle v,f_i\rangle ||^2_A , \forall v\in V,$$ where $b$ is some positive constant and $||.||_A$ is the norm on the C*-algebra $A$.
}\end{defn}

\begin{lem}
If $F$ is a spherical tight $A$-frame in a submodule $V$ of the standard Hilbert $A$-module $\ell^2_A$, then there is an orthonormal basis of $V$, in which $b^{-1/2}F = W_{k,n}U, $ with $U$ is a unitary operator in the $A$-module. 
\end{lem}
\proof  
Indeed let us denote $\mathbf e_1,\dots,\mathbf e_n$ the standard basis of the free $A$-module $A^n$. Then $$\langle b^{-1/2}F\mathbf e_i, b^{-1/2}F\mathbf e_j\rangle = \langle \mathbf e_i, \mathbf e_j\rangle = \delta_{i,j}.$$ Let us denote $V = \langle  b^{-1/2}F\mathbf e_1,\dots, b^{-1/2}F\mathbf e_n \rangle$ the sub$A$-module of the Hilbert $A$-module $B$ and by $V^\perp$ the orthogonal complement
$$V^\perp = \{ v \in V; \langle  v,\mathbf e_i\rangle = \mathbf 0, \forall i=\overline{1,n} \}.$$ Choose orthonormal basis in Hilbert $A$-modules $V$ as follows. 
$$\mathbf f_i= b^{-1/2}F\mathbf e_i, \forall i=\overline{1,n},$$ in $V$  then we have matrix form of the partial isometry 
$$b^{-1/2}F = \left[\begin{array}{cc} b^{-1/2}F|_V & 0\\ 0 & b^{-1/2}F|_{V^\perp}\end{array} \right]=\left[\begin{array}{cc} I_n & 0\\ 0 & 0\end{array} \right]=[I_n|0_{n,k-n}].$$
\endproof

\begin{thm}
Any spherical tight $A$-frame in a finitely generated projective $A$-module $V$ can be realized as a matrix of type
$$F = b^{1/2}W_{k,n}U, \qquad U\in \mathcal O^A(n),$$ where $\mathcal O^A(n)$ is denotes the orthogonal $A$-automorphism group of $A$-module $V$ and $W_{k,n} = [I_n|O_{n,k-n}]$.
\end{thm}
\proof Indeed, 
Every finitely generated projective $A$-module $V$ can be included in some free $A$-module $A^n$, as a direct summand and the latter can be included in the standard Hilbert $A$-module $\ell^2_A$.
From the definition we see that
$$b||v||^2_A = \sum_{i\in I} ||\langle v,f_i\rangle||^2_A .$$ 
This means that 
$b^{-1/2}F $ is a partial isometry in the Hilbert module $A^n$. Following the previous Lemma, we can conclude that in a standard basis we have $b^{-1/2}F = W_{k,n}U , $ where $W_{k,n} = [I_n|O_{n,k-n}]$ and $U $ is an element in the unitary group $U(A)$ with entries from $A$.
\endproof

\section{The structure of the C*-algebra of NC spherical tight $\mathbf A$-frames}

Let us denote by $\mathcal F^A_{k,n}$ the set of all spherical tight $A$-frames in the finitely generated free $A$-module $A^n$.
 
\begin{prop}\label{prop31}
Let $F =[f_1,\dots,f_k]\in \mathcal F^A_{k,n}$ and let $I \subset \{1,\dots,k\}$ be a subset. Let $\mathbf e_1,\dots,\mathbf e_k$, $\mathbf e_i = \left[\begin{array}{c} 0\\ \vdots\\ 1\\ \vdots\\ 0 \end{array}\right]$  be the standard basis of $A^k$, $Q_I: A^k \to A^k$ be the projection onto the subspace $\langle \mathbf e_1,\dots,\mathbf e_k\rangle_A = span_A\{\mathbf e_i | i\in I\}$. Then $F^*F$ is a projector and commutes with the projector $Q_I$ if and only if there is a submodule $V \subseteq A^n$ such that $[f_i]_{i\in I}$ is a tight $A$-frame in sub-$A$-module $V$, and $[f_i]_{i\in I^c}$ is a spherical tight $A$-frame in the sub-$A$-module $V^\perp$. Moreover, in that case (where $F^*F$ commutes with $Q_I$), the cardinality of $I$ is a multiple of $k/d$, where $d= gcd(k,n)$.
\end{prop}
\proof  Without loss of generality, up-to a change of order of $f_i$, which is equal to multiplication by $F$ on the right by a permutation matrix, we may assume that $F=[F_1|F_2]$, i.e.  we may assume that $I=\{1,\dots,p\}$ for some $p\in \{1,\dots, k\}$. Then $F^* = \left[\begin{array}{c} F_1^*\\ F_2^* \end{array}\right]$ and therefore 
$$F^*F = \left[ \begin{array}{cc} F_1^*F_1 & F_1^*F_2\\ F_2^*F_1 & F_2^*F_2 \end{array} \right].$$
In the case $F^*F$ commutes with $Q_I$, we have $F_1^*F_2 = (F_2^*F_1)^* = 0$ and the proof is achieved.
\endproof

\begin{defn}{\rm In the situation of the previous Proposition \ref{prop31}, if $I$ is a proper nonempty subset of $\{1,\dots,k\}$,
we say that the spherical tight $A$-frame $F=[f_1,\dots,f_k]$ is {\it ortho-decomposable}, $[f_i]_{i\in I}$ is  a spherical tight $A$-frame in the sub-$A$-module $V$ and $[f_i]_{i\in I^c}$ is a spherical tight $A$-frame in the sub-$A$-module $V^\perp$. In the opposite case the $A$-frame is called {\it ortho-indecomposable}
}\end{defn}

\begin{thm}
Let $A$ be a C*-algebra, $\sigma$ a partition of the set $\{1,\dots,k\}$, denote $\hat{M}^A_{k,n}$ the set of all $A$-frames $F\in \mathcal F^A_{k,n}$ which are otho-indecomposable, $\hat{M}^A_{k,n}(\sigma)$ the set of all $A$-frames $F\in \mathcal F^A_{k,n}$ such that $\rho_F = \sigma$. 
\begin{enumerate}
\item  $\hat{M}^A_{k,n}$ is a C*-algebra, generated by a ortho-indecomposable system of partial isometries.
\item  Let $d =gcd(k,n)$ and let $k'=k/d$, $n'=n/d$. Let $P(k,k')$ the set of all partitions of the set $\{1,\dots,k\}$ into the subsets whose cardinalities are multiple of $k'$. Then
$\mathcal F^A_{k,n}$ is a direct sum of the algebras $\hat{M}(\sigma)$, for  $\sigma \in \mathcal P(k,k')$,
$$\mathcal F^A_{k,n} = \bigoplus_{\sigma\in \mathcal P(k,k')} \hat{M}^A_{k,n}(\sigma).$$
\item If $\sigma =\{A_1,\dots, A_\ell\} \in \mathcal P(k,k')$ with $|A_i| = m_ik'$, then it is a tensor product of algebras
$$\hat{M}^A_{k,n}(\sigma) \cong \bigotimes_{i=1}^\ell \hat{M}^A_{m_ik',m_in'}$$
\end{enumerate}
\end{thm}
\proof
The theorem is proven in the same way as it was done in the commutative case, see ({loc. cit.},\S 4) for a more detailed analysis.  
The first assertion is exactly deduced from the definition of C*-algebra generated by a system of partial isometries, as the $A$-span of those partial isometries.
\endproof

\section{Application to noncommutative classification spaces}

It is easy to see that $A=\Mat_1(A)$ can be included in $\Mat_n(A)$ and can be regarded as some noncommutative Serre fibration (NCSF), if $A$ admits some NCCW complex structure. The algebra $Mat_\infty(A)$ can then be regarded as some universal NCSF, in the sense that any finite rank NCSF can be obtained as induced one from a universal NCSF
$$\omega_W :\quad \Mat^n_\infty(A) = \lim_{\rightarrow\atop k} \Mat_k^n(A)\to W.$$

\begin{defn}{\rm
Let us consider the space of all infinite matrix of $$\Mat^n_\infty(A) = \lim_{\rightarrow\atop k} \Mat_k^n(A) .$$ Suppose that $\Mat^n_\infty(A) \to W$ is a   $\Mat^n_\infty(A)$-module, then the pushout diagram
$$\CD A \star_{M^n_\infty(A)} W @<<<  W\\
 @AAA   @AAA\\
A @<<< M^n_\infty(A) \endCD$$
 of the two morphisms $M^n_\infty(A) \to A$, defined by $m=[m_{ij}] \mapsto m_{11}$ and the fibration $M^n_\infty(A) \to W$ gives an $A$-module $A \star_{M^n_\infty(A)} W$, which is called {\it NC induced bundle}, what is indeed an induced $A$-module.
}\end{defn}

\begin{thm}
Any rank $n$ projective $A$-module $V$ can be obtained as some induced bundle $A \otimes_{M^n_\infty(A)} W$ from a finitely generated universal projective $\Mat_\infty(A)$-module $\omega_W$.
\end{thm}
\proof
Indeed there is one-to-one correspondence between finitely generated projective $A$-modules and idempotents $e^2 = e^* = e$ in $\Mat_\infty(A)$, in one hand side and one-to-one correspondence between finitely generated projective $M^n_\infty(A)$-modules and idempotents $E^2 = E^* = E$ in $\Mat_\infty(M^n_\infty(A)) \cong M^n_\infty(A) \otimes \Mat_\infty(\mathbf C) \cong \Mat_\infty(A)$ \endproof



Let us denote $A\star B$ the free product of two algebras $A$ and $B$.

\begin{thm}[The Milnor universal bundle] We have a natural noncommutative bundle 
$$\omega_A :\quad \CD   \Mat_\infty(A) \star  \Mat_\infty(A) \star \dots \to  \Mat_\infty(A)\endCD .$$ 
\end{thm}
\proof
It is clear that $\Mat_\infty(A)$ and $\Mat_\infty(A) \star \Mat_\infty(A) \star \dots$ are NCCW. The maps from $\Mat_\infty(A)\star \Mat_\infty(A)\star\dots \to \Mat_\infty(A)$ is defined as the natural product of factors from the free product to the algebra. Up-to homotopy equivalence in the category of NCCW every map is a NCSF.
\endproof

\begin{thm}[Comparison with the free product construction]
There is a natural map from the Milnor universal bundle $\omega_A$ to the universal bundle $\omega_W$, following the commutative diagram

$$\CD
 \Mat_\infty(A) @>>>  W\cong \omega_W = M^n_\infty(A) \star_{\{\Mat_\infty(A)\star \Mat_\infty(A) \star \dots\}} W\\
@A{\omega_A} AA  @AA{\omega_W}A\\
\Mat_\infty(A)\star  \Mat_\infty(A)\star \dots @>>> \Mat^n_\infty(A) \endCD $$
where $W= \mathcal F_n(A^\infty)$ is the universal totaulogical $M^n_\infty(A)$ module.
\end{thm}
\proof
The horizontal arrows are defined following the push-out diagram. They provide the isomorphism between two fibrations in vertical columns.
\endproof

\begin{rem}
The above theorem explains the construction of of J. Cuntz and D. Quillen for classifying NC space as $qA$. We construct some NC analog of the Stiefel construction of classifying spaces.
\end{rem}

\begin{thm}
The K-theory $\K_*(A)$ of $A$ and the $\KK^*(A,\mathbf C)$ are isomorphic.
\end{thm}
\proof
The $\KK$ functor in Cuntz setting is $$\KK(A,B) = [[B,q\Mat_\infty(A)]],$$ where 
$$q\Mat_\infty =\ker\{ \Mat_\infty(A) \star  \Mat_\infty(A) \star \dots \to \Mat_\infty(A)\}$$ and $[[A,B]]$ is the set of all homotopy classes of quasi-isomorphisms from $A$ to $B$.
From the other side, we have $M^n_\infty(A)$ is the classification space
of rank $n$ projective $A$-modules, $$K_*(A) = [\mathbf C,\Mat_\infty(A)].$$ 

We have a natural map see \cite{manuilov-thomsen} from $[B,\Mat_\infty(A)]$ into $[[B,q\Mat_\infty(A)]]$, induced from the natural morphism 
$$q\Mat_\infty(A) = \ker\{\Mat_\infty(A) \star \Mat_\infty(A) \star \dots \to \Mat_\infty(A)\} \to \Mat_\infty(A) .$$ 
And if this morphism becomes isomorphism up to a homotopy, then $K_*(A) \cong KK^*(A,\mathbf C)$.

In one hand we have the KK-theory of the algebras $A$ and $B=\mathbf C$. In  the other hand we have K-theory of $A$ in the sense of the group of stable classes of projective finitely generated $A$-modules over $M^n_\infty(A)$, i.e. the ordinary $K$-homology $\K_*(A)$. The theorem is proved.
\endproof

\begin{defn}{\rm
A {\it spherical tight $A$-frame} in the universal Hilbert $A$-module $\ell^2(A)$ is a set $F=[f_1,\dots,f_k]$ of $k$ vectors in $A^n \subset \ell^2(A)$, such that $k>n$.
}\end{defn}

\begin{cor}
Any spherical tight $A$-frame can be obtained from some spherical tight $A$-frame on the universal Hilbert $A$-module $\ell^2(A)$. 
\end{cor}
\proof In $A$-module $A^n$ $F=[f_1,\dots,f_k]$ is of form 
$F= b^{1/2}W_{k,n}U,$ with $U\in \mathcal O^n(A)$. The latter groups $\mathcal O^A(n)$ is a subgroup in the $\mathcal O^A(\infty)$ and therefore is
also a spherical tight $A$-frame in a $\Mat_\infty(A)$-module $\ell^2(A)$. \endproof

\section*{Acknowledgments}
The author would like to thank Abdus Salam ICTP for an excellent scientific stay and especially to Professor Le Dung Trang for invitation to participate the Commutative and  Noncommutative Geometry Year Program in ICTP.
The author expresses his deep thanks to Professor Claude Schochet for pointing out the work \cite{dykemastrawn}.

\vskip 1cm
{\noindent\sc Institute of Mathematics, Vietnam Academy of Science and Technology, 18 Hoang Quoc Viet Road, Cau Giay District, 10307 Hanoi, Vietnam}

{\noindent\tt Email: dndiep@math.ac.vn}
\end{document}